\documentclass[twoside]{article}
\usepackage{latexsym,amssymb,amsmath}
\usepackage{amsfonts}
\RequirePackage{graphics,epsfig,psfrag}
\def\abs#1{\left \vert #1 \right \vert}

\def\Mod#1{\,(\hbox{\rm mod}\,#1)}

\def\id{{\mathchoice{\rm 1\mskip-4mu l}{\rm 1\mskip-4mu l}{\rm 1\mskip-4.5mu l}{\rm 1\mskip-5mu l}}}%
\def\phi{\varphi}

\def\cC{{\cal C}}
\def\cF{{\cal F}}

\def\TT{{\mathrm T}}

\catcode`\@=11
\def\section{\@startsection{section}{1}{\z@}%
    {-21dd plus-8pt minus-4pt}{10.5dd}
     {\centering\normalsize\bfseries\boldmath}}
\def\subsection{\@startsection{subsection}{2}{\z@}%
    {-21dd plus-8pt minus-4pt}{10.5dd}
     {\normalsize\upshape}}
\def\subsubsection{\@startsection{subsubsection}{3}{\z@}%
    {-13dd plus-8pt minus-4pt}{10.5dd}
     {\normalsize\itshape}}
\def\paragraph{\@startsection{paragraph}{4}{\z@}%
    {-13pt plus-8pt minus-4pt}{\z@}{\normalsize\itshape}}

\catcode`\@=12

\def\pn{\medskip\par\noindent}
\def\bi{\vspace{-2pt}\begin{itemize}\itemsep -2pt plus 1pt minus 1pt}
\def\ei{\end{itemize}\vspace{-4pt}}
\def\bn{\vspace{-2pt}\begin{enumerate}\itemsep -2pt plus 1pt minus 1pt}
\def\en{\end{enumerate}\vspace{-4pt}}
\newcommand{\Pf}{{\em Proof}. }

\newcommand{\EPf}{\hbox{}\hfill$\Box$\vspace{.5cm}}
\def\[#1\]{\begin{eqnarray}#1\end{eqnarray}}
\def\$#1\${\begin{eqnarray*}#1\end{eqnarray*}}

\def\pent#1#2{\lfloor\frac{#1}{#2}\rfloor}

\def\abs#1{\left \vert #1 \right \vert}

\def\frac#1#2{{\textstyle{{#1} \overwithdelims.. {#2}}}}
\def\Frac#1#2{{\displaystyle{{#1} \overwithdelims.. {#2}}}}

\catcode`\@=11
\def\system#1{\left\{\null\,\vcenter{\openup\jot\m@th
\ialign{
\strut\hfil$\displaystyle{##}$&
$\,\displaystyle{{}##}\,$\hfil&&
\strut\hfil$\,\displaystyle{##}$&
$\,\displaystyle{{}##}\,$\hfill
\hfil\crcr#1\crcr}}\right.}

\def\cmatrix#1{\left (
\null\,\vcenter{
\ialign{
\hfil${##}\ $\hfil &
\hfil$\ {##}\ $\hfil&&
\hfil$\ {##}\ $\hfil&
\hfil$\ {##}$\hfil
\crcr#1\crcr}}\right )}

\catcode`\@=12

\catcode`\@=11
\def\@opargbegintheorem#1#2#3{\par\addvspace{6pt plus3pt minus2pt}%
    \def\@tempa{#3}%
    \noindent{\bf #1 #2 \ifx\@tempa\empty\unskip\else\unskip\ (#3).\fi\hskip.5em}\csname#1font\endcsname\ignorespaces
\ignorespaces}

\def\@endtheorem{\par\addvspace{6pt plus3pt minus2pt}}
\def\@begintheorem#1#2#3{\par\addvspace{8pt plus3pt minus2pt}%
              \noindent{\csname#1headfont\endcsname#1\ \ignorespaces#3 #2.}%
              \csname#1font\endcsname\hskip6pt\ignorespaces}
\def\@endtheorem{\par\addvspace{8pt plus3pt minus2pt}\@endparenv}

\catcode`\@=12

\newtheorem{theorem}{Theorem}
\newtheorem{corollary}[theorem]{Corollary}
\newtheorem{lemma}[theorem]{Lemma}
\newtheorem{proposition}[theorem]{Proposition}

\newtheorem{remark}[theorem]{Remark}
\newtheorem{example}[theorem]{Example}

\newtheorem{question}[theorem]{Question}

\date{\today}

\begin{document}
\pagestyle{myheadings}
\markboth{P. -V. Koseleff, D. Pecker}{{\em Fibonacci Knots}}
\title{On Fibonacci knots}
\author{P. -V. Koseleff, D. Pecker}
\maketitle
\begin{abstract}
We show that the Conway polynomials of Fibonacci links
are Fibonacci  polynomials modulo 2. We deduce that, when $ n \not\equiv 0 \Mod 4$ and
$(n,j) \neq (3,3),$
the Fibonacci knot $ \cF _j^{(n)} $ is  not a Lissajous knot.
\end{abstract}
\pn {\bf keywords:}{
Fibonacci polynomials, Fibonacci knots, continued fractions}
\section{Introduction}
Fibonacci knots (or links) were defined by J. C. Turner (\cite{Tu}) as
rational knots with Conway notation $\cC(1,1,\ldots, 1).$
He also considered the generalized Fibonacci knots
$\cF_j ^{(n)} = {\cC} (n,n, \ldots , n),$  where $n$ is a fixed integer and
the sequence $(n,\ldots,n)$ has length $j$.

In this paper we determine the Conway and Alexander polynomials modulo 2
of Fibonacci knots. We show that the Conway polynomial
of a generalized Fibonacci knot is a Fibonacci polynomial modulo $2$.

As an application, we show that if $n \not\equiv 0 \Mod 4$ and $ (n,j) \neq (3,3)$
the Fibonacci knot $ {\cF}_j ^{(n)} $ is not a Lissajous knot.

Our results are obtained by continued fraction expansions.
\section{Conway notation and  Fibonacci knots}\label{cf}
The Conway notation (J. H. Conway, \cite{Co}) is particularly convenient for the
important class of rational (or two-bridge) knots.
The Conway normal form $\cC(a_1, a_2, ..., a_m)$
of a rational knot (or link),
is best explained by the following figure.
\psfrag{a}{\small $a_1$}\psfrag{b}{\small $a_2$}%
\psfrag{c}{\small $a_{m-1}$}\psfrag{d}{\small $a_{m}$}%
\begin{figure}[th]
\begin{center}
{\scalebox{.8}{\includegraphics{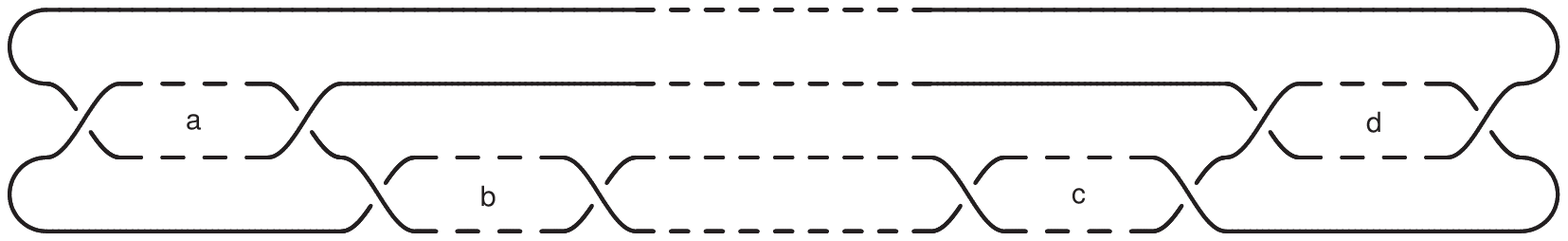}}}\\[30pt]
{\scalebox{.8}{\includegraphics{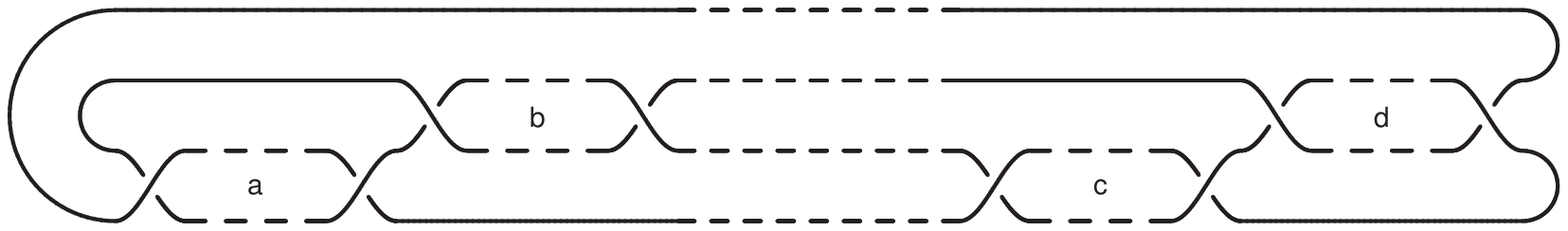}}}
\end{center}
\caption{Conway's normal forms, $m$ odd, $m$ even}
\label{conways3}
\end{figure}
The number of twists is denoted by the integer
$\abs{a_i}$, and the sign of $a_i$ is defined
as follows: if $i$ is odd, then the right twist is positive,
if $i$ is even, then the right twist is negative.
On  Fig. \ref{conways3} the $a_i$ are positive (the $a_1$ first twists are right twists).
\pn
The rational links are classified by their Schubert fractions
\[
\Frac {\alpha}{\beta} =
a_1 + \Frac{1} {a_2 + \Frac {1} {a_3 + \Frac{1} {\cdots +\Frac 1{a_m}}}}=
[ a_1, \ldots , a_m] , \quad \alpha >0.
\]
Two rational links of fractions
$ \Frac {\alpha} {\beta} $ and $  \Frac {\alpha ' }{\beta '} $ are equivalent
if and only if $ \alpha = \alpha' $ and $ \beta' \equiv \beta ^{\pm 1} ( {\rm mod}  \  \alpha).$
The integer $ \alpha$ is the determinant of the link,
it is odd for a knot, and even for a two-component link.
\pn
The following result is a useful consequence of the continued fraction description of
rational links (see \cite{Cr} p. 207).
\begin{theorem}
Any rational link has a Conway normal form $\cC (2a_1, 2a_2, \ldots , 2a_m) .$
\end{theorem}
The Fibonacci knots (or links)  are defined by their Conway notation
$\cF_j = {\cal C} (1,1, \ldots, 1 ), $ where $j$ is the number of crossings.
The Schubert fraction of
$\cF_j $ is $ \Frac{F_{j+1}}{F_j}, $ and its determinant is
the Fibonacci number $F_{j+1}.$
It is the reason why J. C. Turner named these  knots Fibonacci knots.
He also introduced the generalized Fibonacci knots
$\cF _j ^{(n)} = \cC (n,n, \ldots n ),$ where $n$ is a fixed integer.
\begin{figure}[th]
\begin{center}
\begin{tabular}{ccc}
{\scalebox{.6}{\includegraphics{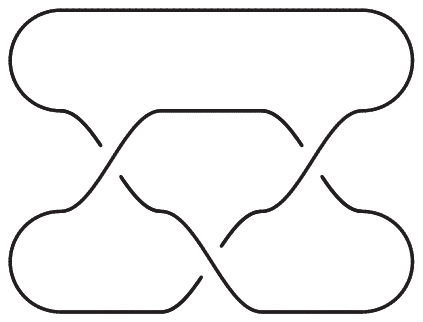}}}&
{\scalebox{.6}{\includegraphics{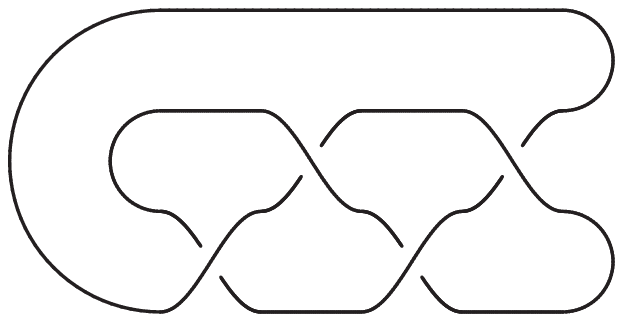}}}&
{\scalebox{.6}{\includegraphics{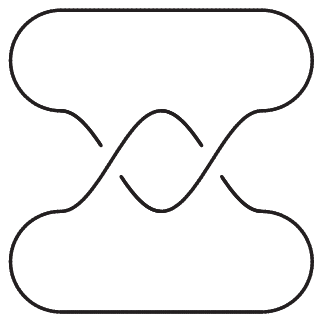}}}\\
$\cF_3^{(1)}$&$\cF_4^{(1)}$&$\cF_1^{(2)}$\\[10pt]
{\scalebox{.6}{\includegraphics{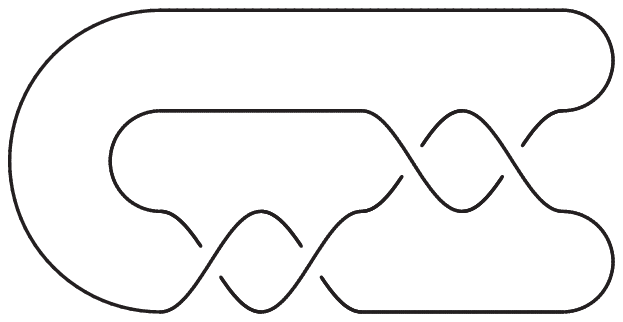}}}&
{\scalebox{.6}{\includegraphics{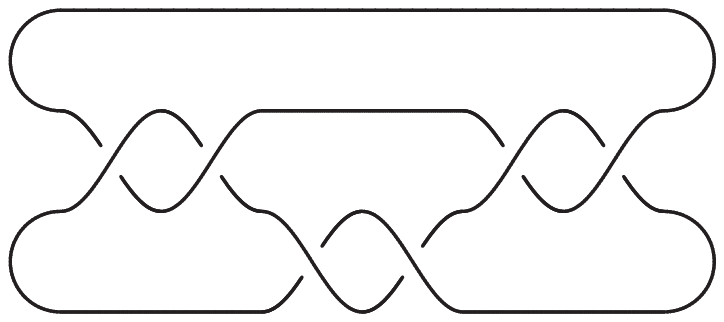}}}&
{\scalebox{.6}{\includegraphics{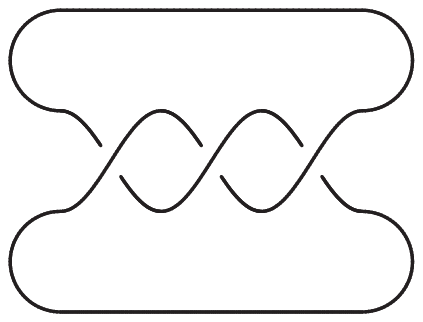}}}\\
$\cF_2^{(2)}$&$\cF_3^{(2)}$&$\cF_1^{(3)}$
\end{tabular}
\pn
\begin{tabular}{cc}
\hfill{\scalebox{.6}{\includegraphics{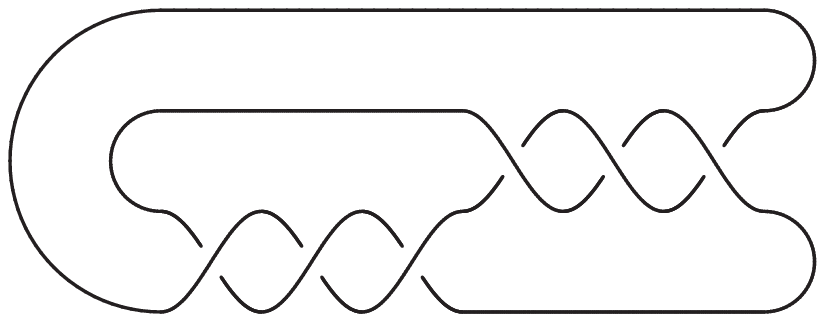}}}&
{\scalebox{.6}{\includegraphics{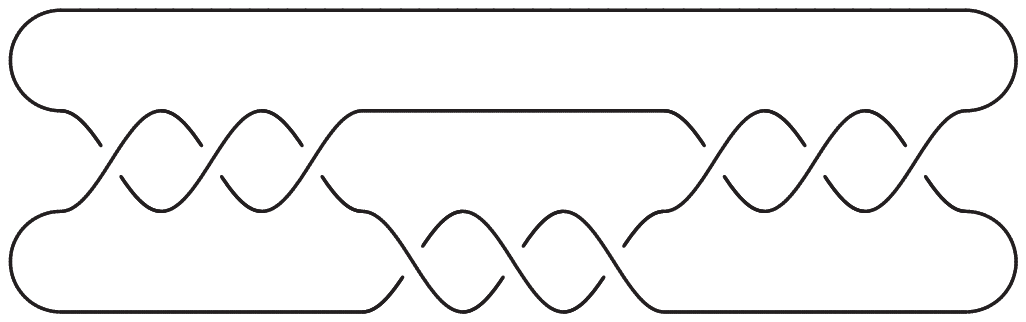}}}\\
$\cF_2^{(3)}$&$\cF_3^{(3)}$
\end{tabular}
\end{center}
\caption{Some Fibonacci knots and links}
\end{figure}
\pn
We first observe
\begin{proposition}\label{fnj}
$\cF_j^{(n)}$ is a knot if and only if
$n \equiv 0 \Mod 2$ and $j \equiv 0 \Mod 2$ or
$n \not \equiv 0 \Mod 2$ and $j \not\equiv 2 \Mod 3$.
\end{proposition}
\Pf
Let us consider the M\"obius transformation $P (z) = [n,z] = n + \Frac{1}{z}$.
It is convenient to consider its matrix notation $P= \cmatrix{n&1\cr1&0}$.
Let $\Frac\alpha\beta = [n,\ldots,n] = P^j (\infty)$, it is also
$$
\cmatrix{\alpha\cr\beta} = P^j \cmatrix{1\cr 0}.
$$
If $n \equiv 1 \Mod 2$ then $P \equiv \cmatrix{1&1\cr1&0}\Mod 2, \ P^2 \equiv \cmatrix{0&1\cr1&1} \Mod 2, \
P^3 \equiv \id \Mod 2$. We deduce that $\alpha \equiv \beta \equiv 1 \Mod 2$ when
$j \equiv 1 \Mod 3$, $\alpha \equiv 0, \, \beta \equiv 1 \Mod 2$ when $j\equiv 2 \Mod 3$ and
$\alpha \equiv 1, \, \beta \equiv 0 \Mod 2$ when $j\equiv 0 \Mod 3$.
The case $n \equiv 0 \Mod 2$ is similar.
\EPf
\section{The Conway and Alexander polynomials}
The Alexander polynomial, discovered in 1928, is one of the most famous invariant of knots.
J. H. Conway discovered an easy way to calculate it.
He introduced the ``Skein relations'' which relate the polynomial of a link $K$
to the polynomials of links obtained by changing one crossing of $K$.
\pn
The following result is a beautiful application of his algorithm.
\begin{theorem}[\cite{Cr}]\label{cr}
Let $K= {\cC} (2a_1, 2a_2, \ldots , 2a_m ) $ be a rational knot (or link).\\
The Conway polynomial of $K$ is
$$ \nabla _K (z)=
\cmatrix{ 1 & 0 } \cmatrix { -a_1 z & 1 \cr 1 & 0}
\cmatrix{ a_2 z & 1 \cr 1 & 0} \cdots
\cmatrix { (-1)^m a_m z & 1 \cr 1 & 0 }
\cmatrix{1 \cr 0 }
$$
The Alexander polynomial of $K$ is
$$ \Delta _K(t)= \nabla _K \bigl( t ^{1/2} - t^{- 1/2}  \bigr).$$
\end{theorem}
Let us consider a simple example.
\begin{example}[The torus links]
The torus link $\TT(2,m)$ has Conway normal form $\cC(m) = \cF_1^{(m)}$.
It is the link of fraction $\Frac m1$ or  $ \Frac m {1-m}. $
We have the continued fraction (of length $m-1$)
$$
\Frac m {1-m} = [ -2,2, \ldots , (-1)^{m-1} \cdot 2 ] .
$$
Hence, the Conway polynomial is
$$
\nabla (z)=
\cmatrix{ 1 & 0} \cmatrix {z & 1 \cr 1 & 0} ^{m -1}\cmatrix { 1 \cr 0 } .
$$
It is well known that
$$
\cmatrix { z & 1 \cr 1 & 0 } ^m =
\cmatrix { f _{m+1} (z) & f_m (z) \cr f_m(z) & f_{m-1} (z) }
$$
where
$f_m(z)$ are the Fibonacci polynomials defined by
$f_0(z)=0, \  f_1(z)=1, \  f_{m+1} (z) = z f_m(z) + f_{m-1}(z) $ (\cite {WP}).
We conclude that the Conway polynomial
of
$\TT(2,m)$ is the Fibonacci polynomial $f_{m} (z)$ (see also \cite {Ka}).
If $m= 2k+1$ ({\em i.e.} $\TT(2,m)$ is a knot) we obtain the Alexander polynomial
$$
\Delta (t)= f_{2k+1} \bigl( t^{1/2} - t^{- 1/2} \bigr)
=
( t^k + t^{-k} ) - ( t^{k-1} + t^{k-1} ) + \cdots + (-1)^{k} .
$$
\end{example}
The recently introduced Lissajous knots
(\cite{BHJS,JP,La,Cr}) are non singular Lissajous space curves.
We will show that in many cases, Fibonacci knots are not Lissajous knots. Let us
first recall the following
\begin{theorem}[\cite{JP,La}]\label{lamm}
If $K$ is  a rational Lissajous knot then
$\Delta_K (t) \equiv 1 \Mod 2.$
\end{theorem}
Consequently, we see that a non trivial torus knot is never a Lissajous knot.
\pn
Moreover, Theorem \ref{cr} provides many examples of knots which are not Lissajous knots.
\begin{corollary}\label{fiber}
Let $b_i \equiv 2 \Mod 4,  m >1.$ The Conway polynomial of
$ {\cC} (b_1,\ldots , b_m) $
is equivalent to $f_ {m+1} (z)  \Mod 2 .$
\end{corollary}
\begin{corollary} If $ n\equiv 2 \Mod 4,$ the modulo $2$ Conway polynomial of
$\cF_j ^{(n)} $ is $ f_{j+1}(z).$
\end{corollary}
Hence these knots are not Lissajous knots by Theorem \ref{lamm}.
\pn
The following result is an immediate consequence of Theorem \ref{cr}.
\begin{corollary}
If $ n \equiv 0 \Mod 4,$ the modulo $2$ Conway polynomial of
$\cF_j ^{(n)} $ is $0$ if $j$ is odd, and $1$ if $j$ is even.
\end{corollary}
It is not known whether the knot $ {\cal F}_2^{(4)}= {\cal C }(4,4) $ is Lissajous or not
(see \cite{BDHZ}).
\section{The modulo 2 Conway polynomial of Fibonacci knots}
We shall now study  the knots $ {\cal F}_j^{(n)} ,$  where $n=2k+1$ is an odd integer.
\begin{lemma}
Let $n=2k+1.$ We have the identities
\[ [n,n,x]= [n+1, \underbrace {-2,2, \ldots, -2,2}_{2k} , -(1+x) ] \label{e1}\]
\[ [n,n,n,z ] = [ n+1, \underbrace{-2,2, \ldots , -2,2}_{2k} , -(n+1) , -z ]. \label{e2} \]
\end{lemma}
\Pf
Let us prove the first formula.
We shall use matrix notations for M\"{o}bius transformations.
Let $G(u)= [-2,2,u]= \Frac {3u+2}{-2u-1}.$
Its matrix is
$ G= \cmatrix{ 3 & 2 \cr -2 & -1 },$ and consequently we get
by induction
\[
G^k= \cmatrix{ 1+2k & 2k \cr -2k & 1-2k } =  \cmatrix{n &  n-1 \cr 1-n & 2-n }.
\label{gk}
\]
Let
\[ M(x)= [n+1, -2,2, \ldots , -2,2,- (1+x) ], \
L(u)= [n+1,u], \  T(x)= -x-1.
\]
The corresponding matrices are
\[ L= \cmatrix{ n+1 & 1 \cr 1 & 0 } , \  T= \cmatrix{ 1 & 1 \cr 0 & -1 }, \
M = L G^k T.
\]
 Consequently
\[
M= \cmatrix{ n+1 & 1 \cr 1 & 0 }
\cmatrix{ n & n-1 \cr 1-n & 2-n }
\cmatrix { 1 & 1 \cr 0 & -1 }
= \cmatrix{n^2+1 & n \cr n & 1 } =
 \cmatrix{ n & 1 \cr 1 & 0 } ^2,
\]
that is $M(x) = [n,n,x]$
which proves  the first identity.
If we substitute $ x= [n, z] $ in Formula (\ref{e1}), we obtain the second identity (\ref{e2}).
\EPf
\begin{corollary}\label{nnn}
Let $n=2k+1.$ We have the continued fractions
$$ [n,n]= [n+1, \underbrace{-2,2, \ldots , -2,2}_{2k}] , \
[n,n,n] = [n+1,  \underbrace{-2,2, \ldots ,  -2,2}_{2k}, -(n+1)].
$$
Let us denote $ [n]_j = [\underbrace{ n, \ldots , n}_j].$
If $ j \not\equiv 1 \Mod 3, $ we get the continued fractions
$$ [n]_{j+3}= [ n+1, \underbrace{-2,2, \ldots , -2,2}_{2k}, -(n+1) , -[n]_j ].$$
\end{corollary}
When $j \equiv 1 \Mod 3$, there is no continued fraction expansion of
$[n]_j$ with even quotients,
by Prop. \ref{fnj}.
In this case, we shall get a continued fraction expansion for
$ \Frac {\alpha }{\beta - \alpha},$ which is another fraction of the same knot.
Let $s$ be the M\"{o}bius transformation defined by $ s(x)= \Frac x {1-x} . $
We have $ s ( \Frac {\alpha}{\beta} ) = \Frac \alpha {\beta - \alpha  } .$
\begin{proposition}\label{[n]}
Let $n=2k+1.$ We have the continued fractions
\[ s ( [n]_1)=s(n)= \Frac n{1-n} = [ \, \underbrace{-2,2, \ldots -2,2 }_{2k} \,] , \  n\neq1,\]
\[
s([n]_{j+3} ) = [\, \underbrace{-2,2, \ldots , -2,2}_{2k}\, ,-( n+1), -(n+1), -s([n]_j)]
\].
\end{proposition}
\Pf
The first formula has already been proved. Let us prove the second formula.
We shall use the M\"{o}bius maps
$ G(x)= [-2,2,x] = \Frac { 3x+2}{-2x-1} , \   Q(x)= [- (n+1), x] , \ R(x)= -x$
corresponding to
$$ G=\cmatrix{ 3 & 2 \cr -2 & -1 } , \
Q = \cmatrix{ - (n+1) & 1 \cr 1 & 0 }, \
R= \cmatrix{ 1 & 0 \cr 0 & -1 } .
$$
Let us define the M\"{o}bius transformation
$H= G^k \cdot Q^2 \cdot R .$
We obtain, using Form. (\ref{gk})
\[
H
=
\cmatrix { n^3 + n^2 + 2n +1 & n^2 +1 \cr -n^3-n & n-n^2-1 }.
\]
Let $S$ be a matrix corresponding to the M\"{o}bius map $s$.
We have
$$
S^{-1} H S =
\cmatrix{1 & 0 \cr 1 & 1 }
H
\cmatrix{1 & 0 \cr -1 & 1 } =
\cmatrix{ n^3 +2n & n^2 +1 \cr n^2 +1 & n }=
\cmatrix{n & 1 \cr 1 & 0 } ^3,
$$
and then
$$
S \cmatrix{n & 1 \cr 1 & 0 } ^3 = HS.
$$
This means that
$ s \bigl( [n,n,n,x] \bigr) = h \circ s (x),$ which proves our formula.
\EPf
\pn
\begin{remark}
By considering the case $n=1$ where
$h(x)= [-2,-2,-x]$, we obtain
the following interesting continued fractions of length $2m$:
$$ \Frac { F_{3m+2} }{F_{3m}}= [2,2, -2,-2, \ldots ,(-1)^{m-1}\cdot  2, (-1)^{m-1}\cdot 2] $$
Using the corollary \ref{nnn} we obtain similarly
$$ \Frac { F_{3m+1}} { F _{3m}} =
[2, -2,-2, 2,2, \ldots , (-1)^{m-1}\cdot 2, (-1)^{m-1}\cdot  2 , (-1)^m\cdot  2 ] $$
of length $2m$ and
$ \Frac {F _{3m+3} }{ F _{3m+2} } = [2,  - \Frac { F _ {3m+2} }{ F_{3m}} ] $
of length $ 2m+1.$
\pn
Of course, these fractions correspond to Fibonacci knots (or
links). They are not Lissajous knots because of corollary \ref{fiber}.
\end{remark}
It is straightforward to calculate the
Conway polynomials of our Fibonacci knots, using Prop. \ref{[n]}.
\begin{theorem}
Let us denote by $\nabla _j^{(n)}(z)$
the modulo 2 Conway polynomial of the Fibonacci link $ {\cal F}_j^{(n)}. $
We have $ \nabla_j^{(n)} (z)= f_N(z)$ where
\[
\left \{
\begin{array}{l}
\hbox{If } n \equiv 1 \Mod 4, \ N= \pent{j+2}3 (n-2) + j+1,
\\
\hbox{If } n \equiv 3 \Mod 4, \ N= \pent{j+2}3 (n+2) -(j+1).
\end{array}
\right .
\]
\end{theorem}
\begin{corollary}
If $n \not\equiv 0 \Mod 4 $ and $( n, j )\neq (3,3),$ the Fibonacci link
$ \cF_j ^{(n)} $ is not a Lissajous knot.
\end{corollary}
It is not known whether the knot $ \cF_3^{(3)} = {\cal C} (3,3,3)$
is a Lissajous knot (\cite{BDHZ}).
\begin{question}
It would be interesting to study the wider classes of knots defined by their
Conway notation
$ {\cal C} ( \pm n, \pm n, \ldots \pm n) .$
\end{question}
If $n=1$ we obtain all the rational knots (\cite{KP3,KP4}).\\
If $n=2$ we obtain the important class of rational fibered knots (see  \cite{Kaw}).
\pn
In general, we obtain knots with fractions $\Frac \alpha \beta $
such that $ ( \alpha, \beta) \equiv (0, \pm 1)  \ {\rm or } \  ( \pm 1, 0 ) \Mod n .$
\pn
{\bf Acknowledgements:} We would like to thank Pr. C. Lamm for having suggested
 this problem to us.

\vfill
\pn
\hrule width 5cm height 2pt
\pn
Pierre-Vincent Koseleff, \\
\'Equipe-project INRIA Salsa \& Universit{\'e} Pierre et Marie Curie (UPMC-Paris 6)\\
e-mail: {\tt koseleff@math.jussieu.fr}
\pn
Daniel Pecker, \\
Universit{\'e} Pierre et Marie Curie (UPMC-Paris 6)\\
e-mail: {\tt pecker@math.jussieu.fr}


\begin{thebibliography}{ZZ99}
%
\bibitem{BDHZ}
A. Boocher, J. Daigle, J. Hoste, W. Zheng,
{\it Sampling Lissajous and Fourier knots},
{\tt arXiv:0707.4210}, 2007
%
\bibitem{BHJS}
M. G. V. Bogle, J. E. Hearst, V. F .R. Jones, L. Stoilov, {\it
Lissajous knots}, Journal of Knot Theory and its Ramifications,
3(2): 121-140 (1994)
%
\bibitem{Co}
J. H. Conway, {\it An enumeration of knots and links, and some of their
algebraic properties,} Computational Problems in Abstract Algebra
(Proc. Conf., Oxford, 1967), 329-358 Pergamon, Oxford (1970)
%
\bibitem{Cr}
P. Cromwell, {\it Knots and links.} Cambridge University Press, 2004
%
\bibitem{JP}
V. F. R. Jones, J. Przytycki, {\it Lissajous  knots and billiard knots,}
Banach Center Publications, 42:145-163 (1998)
%
\bibitem{Ka}
L.H.  Kauffman , {\it On knots},
Annals of Mathematics Studies, 115. Princeton University Press, Princeton, NJ, 1987
%
\bibitem{Kaw}
A. Kawauchi {\it  A survey of knot theory},
 Birkh\"auser Verlag, Basel, 1996
%
\bibitem{KP3}
P. -V. Koseleff, D. Pecker, {\it Chebyshev knots},
{\tt arXiv:0812.1089}, 2008
%
\bibitem{KP4}
P. -V. Koseleff, D. Pecker, {\it Chebyshev diagrams for rational knots}, \\
{\tt arXiv:0906.4083}, 2009
%
\bibitem{La}
C. Lamm, {\it There are infinitely many Lissajous knots,}
Manuscripta Math., 93: 29-37  (1997)
%
\bibitem{Tu}
J.C. Turner, {\it On a class of knots with Fibonacci invariant numbers},
Fibonacci Quart. {\bf 24, 1},  61-66  (1986)
%
\bibitem{WP}
W. A. Webb, E. A. Parberry, {\it Divisibility of Fibonacci polynomials},
Fibonacci Quart. {\bf 7, 5},  457-463  (1969)
\end{thebibliography}
\end{document}